\documentclass[preprint,numbers]{elsarticle}
\usepackage{amssymb,amsfonts,amsmath,mathrsfs}

\newtheorem{thm}{Theorem}

\newtheorem{prop}[thm]{Proposition}
\newtheorem{cor}[thm]{Corollary}
\newdefinition{rmk}[thm]{Remark}
\newdefinition{ex}[thm]{Example}
\newdefinition{df}[thm]{Definition}
\newproof{pf}{Proof}

\def\slantfrac#1#2{\hbox{$\,^#1\!/_#2$}}

\newcommand{\N}{\mathbb N}

\newcommand{\cl}{\textrm{Cl}}

\newcommand{\cf}{\textrm{cf}}
\newcommand{\St}{\textrm{St}}
\newcommand{\HLS}{\mathcal L_{1\!\tiny{\slantfrac{1}{2}}}}
\newcommand{\dom}{\textrm{dom}}
\newcommand{\ran}{\textrm{ran}}

\begin{document}

\begin{frontmatter}
 \title{Star-covering properties: generalized $\Psi$-spaces, countability conditions, reflection}

 \author[lpa]{L.P.~Aiken}
 \ead{laiken1@gmu.edu}

 \address[lpa]{Department of Mathematical Sciences, George Mason University,
    4400 University Drive, Fairfax, VA 22030, USA}

\begin{abstract}

We investigate star-covering properties of $\Psi$-like spaces.  We show star-Lindel\"ofness is reflected by open perfect mappings.  In addition, we offer a new equivalence of CH.

\end{abstract}
\end{frontmatter}

\section{Introduction}

Our research was motivated by the papers \cite{AJW2011,vMTW2007}.  We answer several questions posed in \cite{AJW2011} as well as some closely related questions.

\begin{df}[\cite{vMTW2007}]

For any topological property $\mathcal P$, $X$ has property \textit{star}-$\mathcal P$ if and only if for each open cover $\mathscr U$ of $X$, there exists a subspace $Y\subseteq X$ such that $Y$ has property $\mathcal P$ and $\St(Y,\mathscr U)=X$.

\end{df}

It is well-known, for example, that every topological space is star-discrete.  It is not hard to see that a space is star-countable if and only if it is star-separable (Lemma 2.3 in \cite{AJW2011}).  For an in-depth discussion of a variety of star-$\mathcal P$ properties, see \cite{AJvMTW2011,AJW2011,Matveev1998,vMTW2007}.  We caution readers to check each author's usage of terminology when reading the literature as it varies from author to author.

The theory becomes more interesting when star-covering properties are considered in conjunction with other properties.  In \cite{AJW2011}, the authors investigate, among other things, the relationship between the star-Lindel\"of and star-countable properties, and  pseudocompactness.  Pseudocompactness is particularly interesting in this case as it may be treated as a star-covering property (see e.g. \cite{Matveev1998,vDRRT1991}).

The authors of \cite{AJW2011} showed that the $\Psi$-space construction provides natural examples of spaces with a variety of behavior.  In this article, we answer many of the questions posed in \cite{AJW2011} using $\Psi$-like spaces.  We also leverage a characterization of pseudocompactness in dense subsets of the Cantor Cube (see \cite{Matveev2010,Reznichenko1989}) to give a characterization of star-Lindel\"ofness within the class of dense pseudocompact subspaces of $2^{\mathfrak c}$.

We will start by answering Question 2 of \cite{AJW2011}:  Is a first-countable, star-Lindel\"of space star-countable?  We will use the following well-known proposition.  The proof is offered as a convenience to the reader.
\begin{prop}[Folklore]
\label{CountableClosure}

Suppose $X$ has an uncountable closed discrete subspace $F$ whose points can separated by pairwise disjoint open sets.  Then $X$ is not star-countable.

\end{prop}
\begin{pf}

Choose $F\subseteq X$ as in the hypothesis. For each $x\in F$, let $U_x\subseteq X$ be an open set containing $x$ such that for each $y\in F\setminus\{x\}$, $U_x\cap U_y=\emptyset$.  Then $\mathscr U=\{U_x\,:\,x\in F\}\cup\{X\setminus F\}$ is an open cover for which there is no countable $Y\subseteq X$ such that $\St(Y,\mathscr U)=X$. \qed

\end{pf}

Recall that for a Hausdorff space $X$, the \textit{Alexandroff Duplicate of} $X$, which we denote $AD(X)$, is the topological space whose point-set is  $X\times\{0,1\}$ topologized by the coarsest Hausdorff topology extending $\{U\times\{0,1\}\,:\,U\subseteq X\mbox{ is open}\}$ $\cup\left\{\{\langle x,1\rangle\}\,:\,x\in X\right\}.$

\begin{ex}

Let $X=AD(\mathbb I)\times (\omega+1)\setminus (\mathbb I\times\{0\}\times\{\omega\})$ where $\mathbb I$ denotes the closed unit interval.  It is clear that $X$ is first-countable and Tychonoff.  For any cover $\mathscr U$ of $X$, $\St(Y,\mathscr U)=X$ where $Y=AD(\mathbb I)\times\omega$. Thus $X$ is star-($\sigma$-compact), hence $X$ is star-Lindel\"of.

For $p\in\mathbb I$, let $U_p=\{\langle p,1\rangle\}\times(\omega+1)$.  Then $\{U_p\,:\,p\in\mathbb I\}$ is a pairwise disjoint collection of open sets separating $\{\langle p,1,\omega\rangle\,:\,p\in\mathbb I\}$ which is closed.  By Proposition \ref{CountableClosure}, $X$ is not star-countable.   \qed

\end{ex}
\begin{df}[Iterated Stars]

Suppose $\mathscr A$ is a family of subsets of $X$ and $Y\subseteq X$.  For $n\in\N$, we define $\St^{(n+1)}(Y,\mathscr A)=\St(\St^{(n)}(Y,\mathscr A),\mathscr A)$ where $\St^{(0)}(Y,\mathscr A)=Y$.

\end{df}

The concepts in the following definition are covered in detail in \cite{Matveev1998} using the terminology $n$-star-compact, $n$-star-Lindel\"of, $n \tiny{\slantfrac{1}{2}}$-star-compact, and $n\tiny{\slantfrac{1}{2}}$-star-Lindel\"of.  To avoid confusion, we will adopt the following substitute notation:

\begin{df}

For each $n\in\N^+$, we say a topological space $X$ has property $\mathcal C_n$ $(\mathcal L_n)$ if and only if for every open cover $\mathscr U$ of $X$, the cover $\{\St^{(n)}(x,\mathscr U)\,:\,x\in X\}$ has a finite (countable) subcover.  We will say $X$ has property $\mathcal C_{n \tiny{\slantfrac{1}{2}}}$ $(\mathcal L_{n \tiny{\slantfrac{1}{2}}})$ if and only if for every open cover $\mathscr U$ of $X$, the cover $\{\St^{(n)}(U,\mathscr U)\,:\,U\in\mathscr U\}$ has a finite (countable) subcover.

\end{df}

The following definition will allow us two work with the two common notions of almost disjoint within a single framework (see \cite{DowVaughan2010,Kunen1980}).
\begin{df}[Generalized $\Psi$-space]

Suppose $\lambda\leq\kappa$ are infinite cardinals and $\mathcal E\subseteq [\kappa]^\lambda$ is a maximal almost disjoint\footnote{By almost disjoint, we mean the intersection of distinct elements has size $<\lambda$.} family (m.a.d.f.), where $[\kappa]^\lambda=\{a\subseteq\kappa\,:\,|a|=\lambda\}$.  Let $\Psi(\mathcal E)$ denote the topological space whose point-set is $\kappa\cup\mathcal E$, with the topology generated by isolating each $\alpha\in\kappa$, and the basic open neighborhoods about $E\in\mathcal E$ are all sets of the form $\{E\}\cup (E\setminus F)$ where $F\in [E]^{<\lambda}$.

\end{df}

In all that follows, $\kappa,\lambda$ and $\mathcal E$ are assumed to be as in the above definition, i.e. $\lambda\leq\kappa$ are infinite cardinals and $\mathcal E\subseteq [\kappa]^\lambda$ is a m.a.d.f.  For convenience and to avoid trivial cases, we will also assume that $\mathcal E$ is disjoint from $\kappa$, $\bigcup\mathcal E=\kappa$, and $|\mathcal E|\geq\kappa$.  In sections $2$ and $3$, unless explicitly stated otherwise, $\lambda$ is assumed to be $\aleph_0$ and $\aleph_1$, respectively.

\section{Properties of $\Psi(\mathcal E)$ when $\lambda=\aleph_0$}

Question 1 (Question 3) of \cite{AJW2011} asks if a first-countable feebly compact\footnote{A space is feebly compact if every locally finite collection of open sets is finite.} (pseudocompact Tychonoff) space is star-Lindel\"of.  We answer both questions in the negative.  Moreover, we will show that a Tychonoff pseudocompact space may fail to be star-$\HLS$, which is, in general, weaker than star-Lindel\"ofness.  This will be sharp within the class of $\Psi$-spaces, as we will show our example has property $\mathcal C_2$, and therefore $\mathcal L_2$, the next property in the hierarchy of Lindel\"of-like star-covering properties (see \cite{Matveev1998,vDRRT1991}).

\begin{prop}
\label{LocallyCountableEquivalence}

Suppose $X$ is locally countable.  Then the following are equivalent:

\begin{enumerate}

\item $X$ is star-countable.

\item $X$ is star-Lindel\"of.

\item $X$ has property $\HLS$.

\end{enumerate}

\end{prop}
\begin{pf}

The equivalence is immediate from the fact that in a locally countable space, every Lindel\"of subspace is countable and contained in an open Lindel\"of subspace. \qed

\end{pf}
\begin{prop}[Folklore]

If $X$ is Hausdorff, sequential and $S\subseteq X$, then $|{\rm Cl}(S)|\leq |S|^\omega$.

\end{prop}
\begin{prop}
\label{ClopenCardinality}

Suppose $X$ is Hausdorff, sequential, each $x\in X$ is contained in an open set of cardinality $\leq\mu$, and $S\subseteq X$.  Then there exists a clopen $K\supseteq S$ such that $|K|\leq \mu^\omega|S|^\omega$.

\end{prop}
\begin{pf}

For each $x\in X$, choose an open neighborhood of $x$, $U_x$, of cardinality at most $\mu$.  Let $J_0=S$, and for $\alpha<\omega_1$, define $J_\alpha=\cl\left(\bigcup\{U_x\,:\,\exists\,\beta<\alpha[x\in J_\beta]\}\right)$.  Let $K=\bigcup_{\alpha<\omega_1}J_\alpha$.  By Proposition \ref{CountableClosure}, $|J_\alpha|\leq \mu^\omega|S|^\omega$, thus $|K|\leq \mu^\omega|S|^\omega$.  $K$ is open because for any $x\in J_\alpha$, $U_x\subseteq J_{\alpha+1}\subseteq K$.  If $x\in\cl(K)$, because $X$ is sequential, there exists a countable $L\subseteq K$ such that $x\in\cl(L)$.  Choose $\alpha<\omega_1$ such that $L\subseteq J_\alpha$, and then, by construction, $x\in J_{\alpha+1}\subseteq K$. \qed

\end{pf}

The following proposition can be found in \cite{Bonanzinga1998}.  We offer a different proof.
\begin{prop}
\label{FirstCountableLargeExtent}

If $X$ is Hausdorff first-countable and $e(X)>\mathfrak c$, then $X$ is not star-countable.\footnote{Here $e(X)$ denotes extent of $X$, that is $\sup\{|F|\,:\,F\subseteq X$ is closed and discrete$\}$.}

\end{prop}
\begin{pf}

Fix a closed discrete subset $F\subseteq X$ of cardinality $\mathfrak c^+$.  For each $x\in F$, let $\langle U_{x,n}\,|\,n<\omega\rangle$ enumerate a countable base at $x$ such that $U_{x,n}\cap F=\{x\}$.  For $n<\omega$, let $\mathscr U_n=\{U_{x,n}\,:\,x\in F\}\cup\{X\setminus F\}$.  Suppose that $Y_n\subseteq X$ is countable and $\St(Y_n,\mathscr U_n)=X$.  Since $U_{x,n}$ is the only open set in $\mathscr U_n$ containing $x$, $Y_n$ intersects each $U_{x,n}$. Thus $F\subseteq\cl\left(\bigcup_{n<\omega} Y_n\right)$, contradicting Proposition \ref{CountableClosure}. \qed

\end{pf}

\begin{prop}

The space $\Psi(\mathcal E)$ is first-countable, Tychonoff, pseudocompact and satisfies property $\mathcal C_2$.

\end{prop}
\begin{pf}

It is clear from the definitions that $\Psi(\mathcal E)$ is first-countable, Hausdorff and zero-dimensional, hence $\Psi(\mathcal E)$ is Tychonoff.  In \cite{Matveev1998}, it is shown that every space with a dense relatively countably compact subspace is $\mathcal C_2$, implying property $\mathcal C_{2\!\tiny{\slantfrac{1}{2}}}$, which is equivalent to pseudocompactness for the class of Tychonoff spaces.  We offer a direct proof as a convenience to the reader.

By the maximality of $\mathcal E$, every infinite subset of $\kappa$ has an accumulation point in $\mathcal E$, i.e. $\kappa$ a dense relatively countably compact subspace.  Thus $\Psi(\mathcal E)$ is pseudocompact.  To verify property $\mathcal C_2$, fix an open cover $\mathscr U$ of $\Psi(\mathcal E)$.  Suppose that $i<\omega$, $\alpha_i\in\kappa$ is such that $\alpha_i\notin\bigcup_{j<i}\St(\alpha_j,\mathscr U)$.  By the maximality of $\mathcal E$, there exists $E\in\mathcal E$ such that $E\cap\{\alpha_i\,:\,i<\omega\}$ is infinite.  Choose $V\in\mathscr U$ such that $E\in V$.  By our choice of $E$, there exists $m<n<\omega$ such that $\alpha_m,\alpha_n\in V$, contradicting that $\alpha_n\notin\St(\alpha_m,\mathscr U)$.  Hence, there exists $F\in[\kappa]^{<\omega}$ such that $\kappa\subseteq\St(F,\mathscr U)$ which implies $\Psi(\mathcal E)=\St^{(2)}(F,\mathscr U)$. \qed

\end{pf}

\begin{prop}
\label{StarCountableGlue}

If $\aleph_0\leq\kappa\leq\mathfrak c$, there exists a m.a.d.f. $\mathcal E\subseteq [\kappa]^\omega$ such that $\Psi(\mathcal E)$ is star-countable.

\end{prop}
\begin{pf}

Let $\mathcal C\subseteq [\omega]^\omega$ and $\mathcal D\subseteq [\kappa\setminus\omega]^\omega$ be maximal almost disjoint families such that $|\mathcal C|=|\mathcal D|=\mathfrak c$.  This is possible because $\kappa\leq\mathfrak c$.  Choose a bijection $f:\mathcal C\to\mathcal D$.  Define $\mathcal E=\{C\cup f(C)\,:\,C\in\mathcal C\}$.  If $A\in[\kappa]^\omega$, $|A\cap\omega|=\aleph_0$ or $|A\cap (\kappa\setminus\omega)|=\aleph_0$.  In either case, $A$ has infinite intersection with some element of $\mathcal C\cup\mathcal D$, thus $\mathcal E$ is maximal.  Then $\Psi(\mathcal E)$ is star-countable because $\St(\omega,\mathscr U)=\Psi(\mathcal E)$, for any open cover $\mathscr U$ of $\Psi(\mathcal E)$. \qed

\end{pf}

The following two propositions each provide negative answers to Questions 1 and 3 of \cite{AJW2011}.
\begin{prop}

If $\aleph_0<\kappa\leq\mathfrak c$, then there exists a m.a.d.f. $\mathcal E\subseteq[\kappa]^\omega$ such that $\Psi(\mathcal E)$ is not star-countable.

\end{prop}
\begin{pf}

Choose an uncountable pairwise disjoint family $\mathcal E_0\subseteq [\kappa]^\omega$, and let $\mathcal E$ be a m.a.d.f. extending $\mathcal E_0$. \qed

\end{pf}

\begin{prop}

If $\kappa>\mathfrak c$ then $\Psi(\mathcal E)$ is not star-$\HLS$.

\end{prop}
\begin{pf}

Suppose $\kappa>\mathfrak c$.  By Proposition \ref{LocallyCountableEquivalence}, it suffices to show that there exists an open cover $\mathscr U$ of $\Psi(\mathcal E)$ such that for every star-countable subspace $X\subseteq \Psi(\mathcal E)$, $\St(X,\mathscr U)\neq \Psi(\mathcal E)$.  Using Proposition \ref{ClopenCardinality}, for $\alpha<\mathfrak c^+$, choose pairwise disjoint clopen $S_\alpha\subseteq\Psi(\mathcal E)$ such that $|S_\alpha|=\frak c$.  Define an open cover $$\mathscr U=\{S_\alpha\,:\,\alpha<\mathfrak c^+\}\cup\{\{\alpha\}\,:\,\alpha\in\kappa\}\cup\left\{\{E\}\cup E\,:\,E\in\mathcal E\setminus\bigcup_{\alpha<\mathfrak c^+}S_\alpha\right\}.$$
Suppose $X\subseteq\Psi(\mathcal E)$ is star-countable and $\St(X,\mathscr U)=\Psi(\mathcal E)$.  Each $S_\alpha$ is infinite and clopen, so we can choose $E_\alpha\in S_\alpha\cap\mathcal E$.  As $S_\alpha$ is the only element of $\mathscr U$ containing $E_\alpha$, $Y\cap S_\alpha\neq\emptyset$.  Let $\mathcal P$ be a partition of $\mathfrak c^+$ into intervals of uncountable length such that $|\mathcal P|=\mathfrak c^+$.   By hypothesis, $X$ is star-countable, so by Proposition \ref{CountableClosure}, for each $I\in\mathcal P$, $Y_I=\bigcup\{X\cap S_\alpha\,:\,\alpha\in I\}$ has an accumulation point $A_\alpha\in\mathcal E\cap X$.  Thus, $\{A_\alpha\,:\,\alpha<\mathfrak c^+\}\subseteq\mathcal E$ is a closed discrete subspace of $X$ of cardinality $\mathfrak c^+$, contradicting Proposition \ref{FirstCountableLargeExtent}. \qed

\end{pf}

\begin{rmk}

One can define the following generalization of $\HLS$: $X$ has property $\mu-\HLS$ if and only for every open cover $\mathscr U$ of $X$, $\{\St(U,\mathscr U)\,:\,U\in\mathscr U\}$ has a subcover of cardinality $<\mu$.  Then the above argument shows that if $\kappa>\mu^\omega$, then $\Psi(\mathcal E)$ is not $\mu^+-\HLS$.

\end{rmk}

\section{Properties of $\Psi(\mathcal E)$ when $\lambda=\aleph_1$}

In this section, we discuss two more questions from \cite{AJW2011}, answering one fully and offering a partial solution to the other.  The following proposition and corollary are essentially contained in the analysis of Example 3.3 of \cite{AJW2011}.
\begin{prop}

If $L\subseteq\Psi(\mathcal E)$ is Lindel\"of, then $L\cap\mathcal E$ and $(L\cap\omega_1)\setminus\bigcup(L\cap\mathcal E)$ are both countable.

\end{prop}
\begin{pf}

Suppose otherwise.  Choose $\mathcal A=\{A_\alpha\,:\,\alpha<\omega_1\}\subseteq L\cap\mathcal E$ where the $A_\alpha$ are taken to be distinct.  For $\alpha<\omega_1$, define $U_\alpha=\{A_\alpha\}\cup\left(A_\alpha\setminus\bigcup_{\beta<\alpha}A_\beta\right).$  Note that $U_\alpha$ is open because $\alpha$ is countable and $A_\beta\cap A_\alpha$ is countable when $\beta<\alpha$.  Then $\mathscr U=\{\{\alpha\}\,:\alpha\in L\cap\omega_1\}\cup\{U_\alpha\,:\alpha<\omega_1\}\cup\{\{E\}\cup E\,:\,E\in L\cap (\mathcal E\setminus\mathcal A)\}$ is an open cover of $L$ with no countable subcover because $U_\alpha$ is the only open set in $\mathscr U$ containing $A_\alpha$.  Also, since $L$ is Lindel\"of and $(L\cap\omega_1)\setminus\bigcup(L\cap\mathcal E)$ is closed and discrete, $(L\cap\omega_1)\setminus\bigcup(L\cap\mathcal E)$ must be countable. \qed

\end{pf}

\begin{cor}

Every Lindel\"of subspace of $\Psi(\mathcal E)$ is contained in a Lindel\"of subspace of the form $\bigcup_{i<\omega}\{E_i\}\cup E_i$ where $\{E_i\,:\,i<\omega\}\subseteq\mathcal E$.

\end{cor}

Problem 3.4 of \cite{AJW2011} asked if $\Psi(\mathcal E)$ can be star-Lindel\"of when $\kappa=\lambda=\aleph_1$.  The authors showed that under the additional hypothesis that $|\mathcal E|^\omega=|\mathcal E|$, $\Psi(\mathcal E)$ is not star-Lindel\"of (Example 3.3 \cite{AJW2011}).  The following is a slight sharpening of their result.  The proof here is essentially the same, as the authors only used the fact that each Lindel\"of subspace of $\Psi(\mathcal E)$ is contained in a Lindel\"of subspace of the form described above.

\begin{prop}

Suppose $\Delta\subseteq[\mathcal E]^\omega$ has cardinality $|\mathcal E|$, and $\Delta$ is order dense in the partial order of reverse inclusion, i.e. for each $C\in[\mathcal E]^\omega$ there exists $D\in\Delta$ such that $D\supseteq C$.  Then $\Psi(\mathcal E)$ is not star-Lindel\"of.

\end{prop}
\begin{pf}

Fix $\mathcal E\subseteq [\kappa]^{\omega_1}$ and an order dense $\Delta\subseteq [\mathcal E]^\omega$ such that $|\Delta|=|\mathcal E|$.  Let $\langle D_\alpha\,|\,\alpha<\mu\rangle$ and $\langle E_\alpha\,|\,\alpha<\mu\rangle$ be enumerations of $\Delta$ and $\mathcal E$, respectively.  As in \cite{AJW2011}, we will build a bijection $f:\mathcal E\to\Delta$ such that $E\notin f(E)$, for each $E\in\mathcal E$.

For $\alpha<\mu$, if $\alpha$ is even, choose $\beta$ least such that $E_\beta\notin\dom(f_\delta)$ for $\delta<\alpha$, and choose $\gamma$ least such that $E_\beta\notin D_\gamma$ and $D_\gamma\notin\ran(f_\delta)$ for $\delta<\alpha$.  If $\alpha$ is odd, choose $\gamma$ least such that $D_\gamma\notin\ran(f_\delta)$ for $\delta<\alpha$, and choose $\beta$ least such that $E_\beta\notin D_\gamma$, and $E_\gamma\notin\dom(f_\delta)$ for $\delta<\alpha$.  Then let $f_\alpha=\{\langle E_\beta,<_\gamma\rangle\}\cup\bigcup_{\delta<\alpha} f_\delta$ and set $f=\bigcup f_\alpha$. It is clear from the construction that $f$ is as desired.

Define an open cover $\mathscr U=\{\{\alpha\}\,:\,\alpha<\kappa\}\cup\{U_E\,:\,E\in\mathcal E\}$ where $U_E=\{E\}\cup E\setminus\bigcup f(E)$.  If $L\subseteq\Psi(\mathcal E)$ is Lindel\"of, by the above proposition, there exists $E\in\mathcal E$ such that $M=\bigcup\{\{F\}\cup F\,:\,F\in f(E)\}\supseteq L$.  Thus $\emptyset=U_E\cap M\supseteq U_E\cap L$, thus $E\notin St(L,\mathscr U)$. \qed

\end{pf}

\begin{prop}

For each cardinal $\mu\geq\aleph_1$, $[\mu]^\omega$ has an order-dense set of cardinality at most $\nu=\mu+\sup\{\xi^\omega\,:\,\aleph_1\leq\xi\leq\mu$ is a cardinal of countable cofinality$\}$.

\end{prop}
\begin{pf}

If $\mu=\aleph_1$, $\omega_1$ is dense in $[\omega_1]^\omega$.  If $\cf(\mu)=\aleph_0$, then $[\mu]^\omega$ is dense in itself and is of size $\nu$.  Otherwise, each countable subset of $\mu$ is bounded, and then by inductive hypothesis, for each $\alpha<\mu$, there exists a dense set $D_\alpha\subseteq[\alpha]^\omega$ such that $|D_\alpha|\leq\nu$.  Then $\bigcup_{\alpha<\mu}D_\alpha$ is a dense subset of $[\mu]^\omega$ of cardinality $\leq\nu$. \qed

\end{pf}

\begin{cor}
\label{UncountableNotStarLindelof}

Suppose that for each uncountable $\mu\leq |\mathcal E|$ of countable cofinality, $\mu^\omega\leq|\mathcal E|$.  Then $\Psi(\mathcal E)$ is not star-Lindel\"of.

\end{cor}

Corollary \ref{UncountableNotStarLindelof} implies that if $|\mathcal E|^\omega=|\mathcal E|$ then $\Psi(\mathcal E)$ is not star-Lindel\"of, as shown in \cite{AJW2011}, but can be used to show even more.  For example, it follows from Corollary \ref{UncountableNotStarLindelof} that if $\mathcal E\subseteq [\kappa]^{\omega_1}$ has cardinality $<\aleph_\omega$, then $\Psi(\mathcal E)$ is not star-Lindel\"of.

Unfortunately, as noted by the authors of \cite{AJW2011}, the situation is more complicated than when $\kappa=\lambda=\aleph_0$.  For example, the existence of a m.a.d.f. $\mathcal E\subseteq[\omega_1]^{\omega_1}$ of cardinality $2^{\omega_1}$ is independent of ZFC (see Chapter 8, Exercise B5 of \cite{Kunen1980}), so it is unclear if a `gluing' argument, similar to that of Proposition \ref{StarCountableGlue}, could be generalized.

Problem 3.5 of \cite{AJW2011} asks if it is consistently true that a feebly Lindel\"of\footnote{A space is feebly Lindel\"of if every locally finite collection of open sets is countable.} $P$-space is star-Lindel\"of.  The above corollary provides numerous ZFC examples of feebly Lindel\"of $P$-spaces that are not star-Lindel\"of.

\begin{ex}

If $\kappa^{\omega_1}=\kappa$, then $\Psi(\mathcal E)$ is not star-Lindel\"of.  It is clear that for $\lambda=\aleph_1$, $\Psi(\mathcal E)$ is a $P$-space, and by the maximality of $\mathcal E$, $\Psi(\mathcal E)$ is feebly Lindel\"of.

\end{ex}

\section{Reflection of Star-Covering Properties}

Recall that a topological property $\mathcal P$ is said to be \textit{reflected by a class of mappings} $\mathcal Q$ if $X$ must have property $\mathcal P$ whenever there exists a mapping of class $\mathcal Q$ from $X$ onto a space with property $\mathcal P$.

Question 5 of \cite{AJW2011} asks if star-countability, star-($\sigma$-compactnness), or star-Lindel\"ofness is reflected by perfect, open or closed mappings.  Mapping an uncountable discrete space onto the space with a single point shows that none of these star-covering properties are reflected by open or closed mappings.  The following example shows that star-countability, star-($\sigma$-compactness), and star-Lindel\"ofness are not reflected by perfect mappings.

\begin{ex}

Let $\mathcal E\subseteq [\omega]^\omega$ be a maximal almost disjoint family.  Let $X=\left(\Psi(\mathcal E)\times\{0,1\}\right)\setminus\left(\omega\times\{1\}\right)$.  Let $f:X\to\Psi(\mathcal E)$ be the projection onto the first coordinate.  Then $f$ is perfect, but $X$ is not star-countable since $\mathcal E\times\{1\}$ is an uncountable closed discrete subspace whose points can be separated by pairwise disjoint open sets, contradicting Proposition \ref{CountableClosure}.

\end{ex}
\begin{rmk}

Alternatively, one could apply Theorem 3.7.29 of \cite{Engelking1989}, which states that any property that is hereditary with respect to clopen subspaces and reflected by perfect mappings is hereditary with respect to closed subspaces.\\

The failure of reflection of star-countability was already known.  If $e(X)$ and $c(Y)$\footnote{$c(Y)$ denotes the cellularity of $Y$, $\sup\{|\mathscr C|\,:\,\mathscr C$ is a pairwise disjoint family of open sets$\}$.} are uncountable, then $X\times Y$ is not star-countable.  Then if $Y$ is compact, the projection of $X\times Y$ onto $X$ is open and perfect.  For more details, see Corollary 2.4 of \cite{BonanzingaMatveev2001} and Example 3.3.4 of \cite{vDRRT1991}.

\end{rmk}

\begin{prop}

Suppose $f:X\to Y$ is an open perfect map and $Y$ is star-Lindel\"of.  Then $X$ is star-Lindel\"of.

\end{prop}
\begin{pf}

Fix an open cover $\mathscr U$ of $X$.  For each $y\in Y$, chose a finite $\mathscr U_y\subseteq\mathscr U$ such that $\bigcup\mathscr U_y\supseteq f^{-1}(y)$ and each $U\in\mathscr U_y$ intersects $f^{-1}(y)$.  Define an open cover of $Y$, $\mathscr V=\{V_y\,:\,y\in Y\}$ where $$V_y=Y\setminus f\left[X\setminus\left(f^{-1}\left[\bigcap_{U\in\mathscr U_y}f[U]\right]\cap\bigcup\mathscr U_y\right)\right].$$
Now, if $q\in V_y$, by the definition of $V_y$, $f^{-1}(q)\subseteq\bigcup\mathscr U_y$ and for each $U\in\mathscr U_y$, $U\cap f^{-1}(q)\neq\emptyset$.  Since $Y$ is star-Lindel\"of, we may choose a Lindel\"of subspace $L\subseteq Y$ such that $\St(L,\mathscr V)=Y$.  Let $M=f^{-1}[L]$.  To see that $\St(M,\mathscr U)=X$, choose $x\in X$, $l\in L$ and $y\in Y$ such that $f(x),l\in V_y$.  Choose $U\in\mathscr U_y$ such that $x\in U$, then $f^{-1}(l)$ intersects $U$ so $\St(M,\mathscr U)=X$.

To see that $M$ is Lindel\"of, let $\mathscr U$ be an open cover of $M$.  For each $l\in L$, choose a finite $\mathscr U_l\subseteq\mathscr U$ such that $f^{-1}(l)\subseteq\bigcup\mathscr U_l$.  Let $V_l=Y\setminus f[X\setminus\bigcup\mathscr U_l]$ and then define $\mathscr V=\{V_l\,:\,l\in L\}$.  Choose a countable $S\subseteq L$ such that $\{V_l\,:\,l\in S\}$ covers $L$.  Then $\bigcup\{\mathscr U_l\,:\,l\in S\}$ is a countable subcover of $M$. \qed

\end{pf}

\begin{rmk}

If $X$ is locally compact, then a continuous open surjection with compact fibers is also closed.  Thus, star-Lindel\"ofness is reflected onto locally compact spaces by continuous open mappings with compact fibers.

\end{rmk}

\section{Dense Pseudocompact Subspaces in Dyadic Cubes}

 For $F\in[\kappa]^{<\omega}$ and $f\in 2^F$, let $O_f=\{p\in 2^\kappa\,:\,p\supseteq f\}$.  We will need following three results, which we present without proof.

\begin{prop}[\cite{Arhangelskii1992,EfimovEngelking1965}]

Suppose $\kappa$ is an infinite cardinal and $X\subseteq2^\kappa$ is dense.  Then $X$ is pseudocompact if and only if for each $I\in[\kappa]^\omega$, $\pi_I[X]=2^I$.

\end{prop}

\begin{thm}[\cite{Matveev2010}]

The Continuum Hypothesis is equivalent to the statement:  Every dense pseudocompact subset of $2^{\mathfrak c}$ has a dense Lindel\"of subspace.

\end{thm}

For the proof of sufficiency, Matveev showed that each dense pseudocompact subspace of $2^{\omega_1}$ contains a dense Lindel\"of subspace.  (Proposition 6 of \cite{Matveev2010})  An obvious corollary is that each dense pseudocompact subspace of $2^{\omega_1}$ is star-Lindel\"of.

\begin{ex}[Reznichenko's Example \cite{Reznichenko1989}]

Let $\mathscr P$ be a partition of $\mathfrak c$ into sets of cardinality $\mathfrak c$ such that $|\mathscr P|=\frak c$.  Let $\langle s_\alpha\,|\,\alpha<\mathfrak c\rangle$ and $\langle P_\alpha\,:\,\alpha<\mathfrak c\rangle$ be enumerations of $\bigcup\{2^S\,:\,S\in[\mathfrak c]^\omega\}$ and $\mathscr P$, respectively.  Define $x_\alpha:\mathfrak c\to\{0,1\}$ by $x_\alpha(\beta)=s_\alpha(\beta)$ if $\beta\in\dom(s_\alpha)$ and $x_\alpha(\beta)=\chi_\alpha(\beta)$ otherwise, where $\chi_\alpha$ denotes the characteristic function of $P_\alpha$.  Then $X=\{x_\alpha\,:\,\alpha<\mathfrak c\}$ is a dense pseudocompact subset of $2^{\mathfrak c}$ such that for each $Y\subseteq X$ of cardinality $<\mathfrak c$, $Y$ is closed and discrete.

\end{ex}

\begin{prop}

The following statements are equivalent:

\begin{enumerate}

\item CH.

\item Every dense pseudocompact subspace of $2^{\mathfrak c}$ is star-Lindel\"of.

\item Reznichenko's Example is star-Lindel\"of.

\end{enumerate}

\end{prop}
\begin{pf}

$(1)\Longrightarrow (2)\Longrightarrow (3)$ follows from the above results.  For $(3)\Longrightarrow (1)$, let $X$ denote Reznichenko's Example. Suppose $Y\subseteq X$ is countable.  Let $\Gamma=\bigcup_{x_\alpha\in Y}\dom(s_\alpha)$ and choose $\delta<\mathfrak c$ such that: $\dom(s_\delta)\cap\Gamma=P_\delta\cap\Gamma=\emptyset$, $s_\delta(\alpha)=0$ for each $\alpha\in\dom(s_\delta)\cap\bigcup_{s_\alpha\in Y}P_\alpha$, and $s_\delta(\alpha)=1$ otherwise.  This choice of $\delta$ is possible because there are $\mathfrak c$-many such functions.  For $\alpha<\mathfrak c$, let $f_\alpha=\{\langle\alpha,1\rangle\}$ and define $\mathscr U=\{O_{f_\alpha}\,:\,\alpha<\mathfrak c\}$.  Fix $x_\alpha\in Y$ and $\beta<\mathfrak c$.  If $\beta\in\dom(s_\delta)$ then $\beta\notin\Gamma$, so $x_\alpha(\beta)=\chi_\alpha(\beta)$, and if $x_\delta(\beta)=1$, then by construction $x_\alpha(\beta)=0$.  If $\beta\in P_\delta\setminus\dom(s_\delta)$, then $x_\alpha(\beta)=0$ because $(P_\alpha\cup\dom(s_\alpha))\cap P_\delta=\emptyset$.  Thus $x_\delta\notin\St(Y,\mathscr U)$.  Assume $\mathfrak c>\aleph_1$.  If $Y\subseteq X$ is such that $\St(Y,\mathscr U)=X$, then $Y$ is uncountable.  Then $Y$ contains a subspace $Z$ of cardinality $\aleph_1<\mathfrak c$.  Thus $Z$ is closed and discrete, and so $Y$ is not Lindel\"of. \qed

\end{pf}
It should be noted that there are dense pseudocompact subsets of $2^{\mathfrak c^+}$ that are not star-Lindel\"of (Remark 3 of \cite{Matveev2010}).

\begin{rmk}

We may not replace star-Lindel\"of with star-countable in the above Proposition.  The proof shows that, irrespective of CH, Reznichenko's Example is not star-countable.  Alternatively, in \cite{Pavlov2005}, it is shown that Reznichenko's Example is meta-Lindel\"of.\footnote{A space is meta-Lindel\"of if every open cover has a point-countable open refinement.}  It is well-known that a meta-Lindel\"of, star-countable space is Lindel\"of, and it is not hard to see that Reznichenko's Example is not Lindel\"of.  The following is a more elementary example of a dense pseudocompact subset of $2^{\mathfrak c}$ that is not star-countable.

\end{rmk}
\begin{ex}

Let $X\subseteq 2^\mathfrak c$ be the $\Sigma$-product with its center, $\mathbf{0}$, removed, where $\mathbf 0$ denotes the constant function taking value $0$.  It is clear that $X$ is dense and pseudocompact.  Let $\mathscr U=\{O_{f_\alpha}\,:\,\alpha<\mathfrak c\}$ be as above.  Then if $Y\subseteq X$ is countable, there exists $\alpha<\mathfrak c$ such that for each $\alpha<\beta<\mathfrak c$ and $p\in Y$, $p(\beta)=0$.  It follows easily that $\St(Y,\mathscr U)\neq X$.  Note that $\mathfrak c$ can be replaced with any uncountable cardinal.

\end{ex}

\section{Problems Remaining Open}

\begin{enumerate}

\item Is $\Psi(\mathcal E)$ star-Lindel\"of when $\kappa=\lambda=\aleph_1$?  More generally, for which $\kappa$, $\lambda$ and $\mu$ is $\Psi(\mathcal E)$ star-$\mu$-Lindel\"of?

\item Is a normal feebly Lindel\"of space star-Lindel\"of?  Is a normal star-Lindel\"of space star-countable?

\end{enumerate}

Question 2 is from \cite{AJW2011}, and the methods developed here appear to have little bearing on the question as the $\Psi(\mathcal E)$ spaces are not normal.\\

The author would like to express his gratitude to Mikhail Matveev and Ronnie Levy for sharing their time and expertise throughout the preparation of this article.

\bibliographystyle{elsart-num-sort}
\bibliography{bibfile}

\end{document}